\documentclass[12pt]{amsart}
\usepackage{amsmath,amsthm,latexsym,amscd,amsbsy,amssymb,amsfonts,fleqn,leqno,
euscript, pb-diagram}

\newtheorem{theorem}{Theorem}[section]

\newtheorem{cor}[theorem]{Corollary}

\theoremstyle{definition}

\theoremstyle{remark}

\numberwithin{equation}{section}

\begin{document}

\title[Another approach to parametric Bing and Krasinkiewicz  maps]
{Another approach to parametric Bing and Krasinkiewicz maps}

\author[V. Valov]{Vesko  Valov}
\address{Department of Computer Science and Mathematics, Nipissing University,
100 College Drive, P.O. Box 5002, North Bay, ON, P1B 8L7, Canada}
\email{veskov@nipissingu.ca}
\thanks{The author was partially supported by NSERC Grant 261914-08.}



\date{}


\begin{abstract}
Using a factorization theorem due to Pasynkov \cite{bp3} we provide
a short proof of the existence and density of parametric Bing and
Krasinkiewicz maps. In particular, the following corollary is
established: Let $f\colon X\to Y$ be a surjective map between
paracompact spaces such that all fibers $f^{-1}(y)$, $y\in Y$, are
compact and there exists a map $g\colon X\to\mathbb I^{\aleph_0}$
embedding each $f^{-1}(y)$ into $\mathbb I^{\aleph_0}$. Then for
every $n\geq 1$ the space $C^*(X,\mathbb R^n)$ of all bounded
continuous functions with the uniform convergence topology contains
a dense set of maps $g$ such that any restriction $g|f^{-1}(y)$,
$y\in Y$, is a Bing and Krasinkiewicz map.
\end{abstract}

\maketitle\markboth{}{Bing, Krasinkiewicz and regularly branched
maps}



\section{Introduction}
All spaces in the paper are assumed to be paracompact and all maps
continuous. All maps from $X$ to $M$ are denoted by $C(X,M)$.
Usually, $C(X,M)$ will carry either the uniform convergence topology
or the source limitation topology. When $X$ is compact, these two
topologies coincide. Unless stated otherwise, a space (resp.,
compactum) means a metrizable space (resp., compactum).

In this paper we provide another approach to prove results
concerning parametric Bing and  Krasinkiewicz maps. The approach is
based on Pasynkov's technique developed in \cite{bp} and \cite{bp3}.

Bing maps and Krasinkiewicz maps have been extensively studied
recent years (see \cite{km}, \cite{kr}, \cite{kr1}, \cite{lev},
\cite{ll}, \cite{m1}, \cite{m2}, \cite{m3}, \cite{mv}, \cite{vv},
\cite{vv1}). Recall that a map $f$ between compact spaces is said to
be a {\em Bing map} \cite{lev} provided all fibers of $f$ are Bing
spaces. Here, a compactum is a {\em Bing space} if each of its
subcontinua is indecomposable. Following Krasinkiewicz \cite{kr}, we
say that a space $M$ is a {\em free space} if for any compactum $X$
the function space $C(X,M)$ contains a dense subset consisting of
Bing maps. The class of free spaces is quite large, it contains all
$n$-dimensional manifolds ($n\geq 1$) \cite{kr}, the unit interval
\cite{lev}, all locally finite polyhedra \cite{st}, all manifolds
modeled on the Menger cube $M_{2n+1}^n$ or the N\"{o}beling space
$N_{2n+1}^n$ \cite{st}, as well as all $1$-dimensional locally
connected continua \cite{st}.

Next theorem follows from the proof of \cite[Theorem 1.2]{vv} where
the special case with $X$ and $Y$ being metrizable was established.

\begin{theorem}
Let $M$ be a free $ANR$-space and $f\colon X\to Y$ be a perfect map
with $W(f)\leq\aleph_0$, where $X$ and $Y$ are paracompact. Then the
maps $g\in C(X,M)$ such that all restrictions $g|f^{-1}(y)$, $y\in
Y$, are Bing maps form a dense set $B\subset C(X,M)$ with respect to
the source limitation topology. Moreover, $B$ is $G_\delta$ provided
$Y$ is first countable.
\end{theorem}

Here, $W(f)\leq\aleph_0$ (see \cite{bp}) means that there exists a
map $g\colon X\to\mathbb I^{\aleph_0}$ such that $f\triangle g\colon
X\to Y\times\mathbb I^{\aleph_0}$ is an embedding. For example
\cite[Proposition 9.1]{bp}, $W(f)\leq\aleph_0$ for any closed map
$f\colon X\to Y$ such that $X$ is a metrizable space and every fiber
$f^{-1}(y)$, $y\in Y$, is separable.

Although,  the arguments from \cite{vv} don't work when the map $f$
in Theorem 1.1 is not perfect or the space $M$ is not $ANR$, we have
the following result:

\begin{theorem}
Let $X$ and $Y$ be paracompact spaces and $f\colon X\to Y$ be a map
with compact fibers and $W(f)\leq\aleph_0$. Then for every compact
free space $M$ the space $C(X,M)$ equipped with the uniform
convergence topology contains a dense subset of maps $g$ such that
all restrictions $g|f^{-1}(y)$, $y\in Y$, are Bing maps.
\end{theorem}

The second type of results concern Krasinkiewicz maps. A space $M$
is said to be a {\it Krasinkiewicz space} \cite{mv} if for any
compactum $X$ the function space $C(X,M)$ contains a dense subset of
Krasinkiewicz maps. Here, a map $g\colon X\to M$, where $X$ is
compact, is said to be Krasinkiewicz \cite{ll} if every continuum in
$X$ is either contained in a fiber of $g$ or contains a component of
a fiber of $g$. The class of Krasinkiewicz spaces contains all
Euclidean manifolds and manifolds modeled on Menger or N\"{o}beling
spaces, all polyhedra (not necessarily compact), as well as all
cones with compact bases (see \cite{kr1}, \cite{ll}, \cite{m1},
\cite{m3}, \cite{mv}).

\begin{theorem}
Let $f\colon X\to Y$ be a map with compact fibers and
$W(f)\leq\aleph_0$, where $X$ and $Y$ are paracompact spaces. If $M$
is a compact Krasinkiewicz space, then $C(X,M)$ equipped with the
uniform convergence topology contains a dense subset of maps $g$
such that all restrictions $g|f^{-1}(y)$, $y\in Y$, are
Krasinkiewicz maps.
\end{theorem}

\begin{cor}
Let $f\colon X\to Y$ be a map with compact fibers such that
$W(f)\leq\aleph_0$, where $X$ and $Y$ are paracompact spaces. Then
for every $n\geq 1$ the space $C^*(X,\mathbb R^n)$ of all bounded
continuous functions with the uniform convergence topology contains
a dense set of maps $g$ such that any $g|f^{-1}(y)$, $y\in Y$, is a
Bing and Krasinkiewicz map.
\end{cor}

Theorem 1.3 was established in \cite[Theorem 1.1]{vv} for an
arbitrary Krasinkiewicz $ANR$-space $M$ in the case $X,Y$ are
metrizable, $f$ is perfect and $C(X,M)$ is equipped with the source
limitation topology. Let us note that the proof of \cite[Theorem
1.1]{vv} provides the following result: Let $f\colon X\to Y$ be a
perfect map between paracompact spaces with $W(f)\leq\aleph_0$, and
let $M$ be a Krasinkiewicz $ANR$-space. Then the maps $g\in C(X,M)$
such that all $g|f^{-1}(y)$, $y\in Y$, are Krasinkiewicz maps form a
dense subset of $C(X,M)$ with respect to the source limitation
topology. Moreover, this set is $G_\delta$ if $Y$ is first
countable.

\textbf{Remark.} The requirement in Theorems 1.2 - 1.3 $f$ to have
compact fibers is necessary because of the definition of Bing and
Krasinkiewicz maps. If we define a Bing space to be a space such
that any its subcontinuum (a connected compactum) is indecomposable,
and a Bing map to be a map whose fibers are Bing spaces, then
Theorem 1.2 remains valid for any map $f$ with $W(f)\leq\aleph_0$.
Then same remark is true for Theorem 1.3 if by a Krasinkiewicz map
we mean any map $g\colon X\to M$, where $X$ is not necessary
compact, such that every continuum in $X$ is either contained in a
fiber of $g$ or contains a component of a fiber of $g$.



\section{Bing and Krasinkiewicz maps}

This section contains the proof of Theorem 1.2, Theorem 1.3 and
Corollary 1.4.

\medskip
\textit{Proof of Theorem $1.2$}. Since $W(f)\leq\aleph_0$, there
exists a map $\lambda\colon X\to\mathbb I^{\aleph_0}$ such that
$f\triangle\lambda\colon X\to Y\times\mathbb I^{\aleph_0}$ is an
embedding. We also fix a map $g_0\colon X\to M$ and a number
$\epsilon>0$. We are going to find a map $g\in C(X,M)$ such that $g$
is $\epsilon$-close to $g_0$ and all restrictions $g|f^{-1}(y)$,
$y\in Y$, are Bing maps. To this end, let
$\overline{\lambda}\colon\beta X\to\mathbb I^{\aleph_0}$ and
$\overline{g}_0\colon\beta X\to M$ be the Stone-Cech extensions of
the maps $\lambda$ and $g_0$, respectively. Then
$\overline{\lambda}\triangle\overline{g}_0\in C(\beta X,\mathbb
I^{\aleph_0}\times M)$. We consider also the constant maps
$\eta_1\colon\mathbb I^{\aleph_0}\times M\to Pt$ and
$\eta_2\colon\beta Y\to Pt$, where $Pt$ is the one-point space.
According to Pasynkov's factorization theorem \cite[Theorem
13]{bp3}, there exist metrizable compacta $K$, $T$ and  maps
$f^{*}\colon K\to T$, $\xi_1\colon\beta X\to K$, $\xi_2\colon K\to
\mathbb I^{\aleph_0}\times M$ and $\eta^{*}\colon\beta Y\to T$ such
that:
\begin{itemize}
\item $\eta^{*}\circ\beta f=f^{*}\circ\xi_1$;
\item $\xi_2\circ\xi_1=\overline{\lambda}\triangle\overline{g}_0$;
\end{itemize}
If $p\colon\mathbb I^{\aleph_0}\times M\to\mathbb I^{\aleph_0}$ and
$q\colon\mathbb I^{\aleph_0}\times M\to M$ denote the corresponding
projections, we have
$$p\circ\xi_2\circ\xi_1=\overline{\lambda}\hbox{~}\mbox{and}\hbox{~}q\circ\xi_2\circ\xi_1=\overline{g}_0.\leqno{(1)}$$
Since $M$ is a free space, there exists a Bing map $\phi\colon K\to
M$ such that $\phi$ is $\epsilon$-close to $q\circ\xi_2$. Then the
map $\overline{g}=\phi\circ\xi_1$ is $\epsilon$-close to
$\overline{g}_0$. Hence, the maps $g=\overline{g}|X$ and $g_0$ are
also $\epsilon$-close. According to $(1)$, we have
$\lambda=\big(p\circ\xi_2\circ\xi_1\big)|X$. This implies that
$\xi_1$ embeds each fiber $f^{-1}(y)$, $y\in Y$, into $K$ (recall
that $\lambda$ embeds the fibers $f^{-1}(y)$ into $\mathbb
I^{\aleph_0}$). Consequently, $f^{-1}(y)\cap g^{-1}(z)$ is
homeomorphic to a subset of $\phi^{-1}(z)$ for all $z\in M$. Since
the fibers $\phi^{-1}(z)$, $z\in M$, are Bing spaces, so are the
spaces $f^{-1}(y)\cap g^{-1}(z)$. Therefore, each restriction
$g|f^{-1}(y)$, $y\in Y$, is a Bing map. \hfill$\square$

\medskip
\textit{Proof of Theorem $1.3$}. We follow the notations and the
arguments from the proof of Theorem 1.2. The only difference now is
that $M$ is a Krasinkiewicz space. So, there exists a Krasinkiewicz
map $\phi\colon K\to M$ such that $\phi$ is $\epsilon$-close to
$q\circ\xi_2$. Then map $g\in C(X,M)$ is $\epsilon$-close to $g_0$
and $f^{-1}(y)\cap g^{-1}(z)$ is homeomorphic to a compact subset of
$\phi^{-1}(z)$ for all $y\in Y$ and $z\in M$. Since $\phi$ is a
Krasinkiewicz map on $K$, for every continuum $C\subset f^{-1}(y)$
we have either $C$ is contained in $\phi^{-1}(z)$ or contains a
component of $\phi^{-1}(z)$ for some $z\in M$. This implies that $C$
is either contained in $f^{-1}(y)\cap g^{-1}(z)$ or contains a
component of $f^{-1}(y)\cap g^{-1}(z)$ for some $z\in M$. Therefore,
any restriction $g|f^{-1}(y)$ is a Krasinkiewicz map.
\hfill$\square$

\medskip
\textit{Proof of Corollary $1.4$}. Let $g_0\colon X\to\mathbb R^n$
and let $\overline{g}_0\colon\beta X\to\mathbb R^n$ be its
Stone-Cech extension. Proceeding as above, we need the following
fact (see \cite{lev} and \cite{ll}, or \cite{vv}): if $M$ is a free
Krasinkiewicz space, then the maps $g\in C(K,M)$ which are both Bing
and Krasinkiewicz form a dense $G_\delta$-subset of $C(K,M)$. Since
$\mathbb R^n$ is both a free space and a Krasinkiewicz space, we can
choose a Bing and Krasinkiewicz map $\phi\in C(K,\mathbb R^n)$ which
is $\epsilon$-close to $\overline{g}_0$. Then any restriction  map
$g|f^{-1}(y)$ is also Bing and Krasinkiewicz. \hfill$\square$



\bibliographystyle{amsplain}

\end{document}